\documentstyle[epsf]{amsart}
\newcounter{pictype}
\def\epsfsize#1#2{\ifnum\value{pictype}=0 0.7#1\else\epsfxsize\fi}

\newcommand{\ppic}[1]{\setcounter{pictype}{0}
\mbox{$\vcenter{\epsffile{#1}}$}
\setcounter{pictype}{1}}

\newcommand{\plainpic}[1]{\epsffile{#1}}

\newcommand{\cpic}[1]{\mbox{$\vcenter{\plainpic{#1}}$}}

\newcommand{\ypic}[2]{\mbox{$\vcenter{\epsfysize=#1\baselineskip\epsffile{#2}}$}}
\newcommand{\xpic}[2]{\mbox{$\vcenter{\epsfxsize=#1\textwidth\epsffile{#2}}$}}

\newtheorem{theorem}{Theorem}
\newtheorem{prop}{Proposition}
\newtheorem{cor}{Corollary}
\newtheorem{deff}{Definition}
\newtheorem{lemma}{Lemma}

\newtheorem{ex}{Example}

\newfont{\bb}{msbm10}

\newcommand{\rtva}{{\mbox{\bb R}}^2}

\newcommand{\sgn}{\mbox{sign}}

\begin{document}
\title{Invariants of knot diagrams and relations among 
Reidemeister moves} 
\author{Olof-Petter \"Ostlund}
\address{Department of Mathematics, Uppsala University, Box 480, S-751
06 Uppsala, Sweden}
\email{olleo@@math.uu.se}
\begin{abstract}
In this paper a classification of Reidemeister moves, which is the most
refined, is introduced. 
In particular, this classification
distinguishes some $\Omega_3$-moves that only differ in how the 
three strands that are involved in the move are ordered on the knot.

To transform knot diagrams of isotopic knots into each other one must
in general use $\Omega_3$-moves of at least two different classes.
To show this, knot diagram invariants that jump only under 
$\Omega_3$-moves are introduced.

Knot diagrams of isotopic knots can be connected by a sequence of
Reidemeister moves of only six, out of the total of 24, classes.
This result can be applied in knot theory to simplify proofs of invariance
of diagrammatical knot invariants. In particular, a 
criterion for a function on Gauss diagrams to define a knot invariant
is presented.
\end{abstract}
\maketitle

\section{Knot diagrams and Reidemeister moves.}\label{sec:defs}
\subsection{Knot diagrams.}
\begin{figure}\caption{The figure eight knot diagram.}\label{fig:figeight}
\ypic{5}{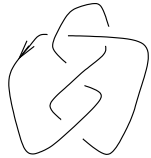}
\end{figure}

A {\em knot diagram} is a picture of an oriented smooth knot like the {\em 
figure eight knot diagram} in Figure~\ref{fig:figeight}. Formally, this is the
image of an immersion of $S^1$ in $\rtva$, with transversal double
points and no points of higher multiplicity, which has been decorated
so that we can distinguish:
\begin{description}
\item[a]{An orientation of the strand, and}
\item[b]{an overpassing and an underpassing strand at each double
point.}
\end{description}
Of course, the manner of decoration, such as where an arrow that indicates the 
orientation is placed, does not matter: 
Two knot diagrams are the same if they are made from
the same image, have the same orientation, and have the same
over-undercrossing information at each crossing point. 

Two knot diagrams are {\em
equivalent} if there is an orientation-preserving diffeomorphism of
the plane that takes one diagram to the other diagram.

\subsection{Reidemeister moves.}
A {\em Reidemeister move} is a transformation of a knot diagram which
looks like one of the transformations in Table~\ref{table:moves} 
inside a disk in $\rtva$
(up to diffeomorphisms)
and leaves the knot diagram 
unchanged outside the disk. 
More formally, we say that the knot diagrams $k$ and $l$
are related by a 
Reidemeister move if they are the same outside $C\subset\rtva$ 
and there are diffeomorphisms $g,h$ from $C$ to a disk $D$ such that 
\begin{description}
\item[a]{$h^{-1}\circ g=Id$ on the boundary of $C$, and}
\item[b]{$g(k\cap C)$ and $h(l\cap C)$ is one of the unordered pairs of
diagram-parts in Table~\ref{table:moves}.}
\end{description}
We shall call $C$ the {\em changing disk} of the Reidemeister move.

\begin{table}\caption{Reidemeister moves.}\label{table:moves}
$$
\begin{array}{cc}
  \begin{array}{ccccc}
    \cpic{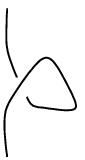}&\begin{array}{c}\rightleftharpoons\\
    \Omega_1\end{array}&
    \cpic{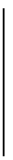}&\begin{array}{c}\rightleftharpoons\\
    \Omega_1\end{array}&\cpic{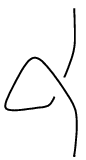}
  \end{array}&
  \begin{array}{ccc}
    \cpic{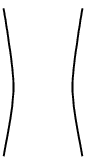}&\begin{array}{c}\rightleftharpoons\\
    \Omega_2\end{array}&\cpic{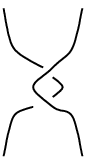}
  \end{array}\\

  \begin{array}{ccc}
    \cpic{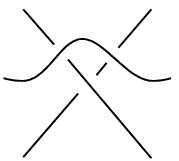} &\begin{array}{c}\rightleftharpoons\\
    \Omega_3\end{array} &\cpic{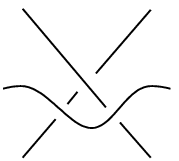}
  \end{array}&
  \begin{array}{ccc}
    \cpic{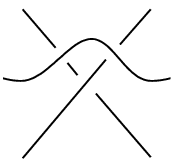} &\begin{array}{c}\rightleftharpoons\\
    \Omega_3\end{array} &\cpic{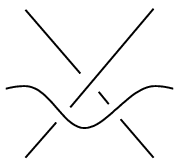} 
  \end{array}
\end{array}
$$
\end{table}

Reidemeister moves were introduced in the setting of PL-knots by
K. Reidemeister~\cite{reidemeister} in 1927. Reidemeister proved the
PL version of the following well-known fact:
{\em Two knot diagrams describe isotopic knots if and only if they are
connected by a sequence of Reidemeister moves and
orientation-preserving diffeomorphisms of the plane.}

\subsection{Classifying Reidemeister moves.}
Reidemeister moves are usually sorted into three classes $\Omega_1$,
$\Omega_2$ and
$\Omega_3$, as in Table~\ref{table:moves}. That is, a move is of
class
$\Omega_i$ if the changing disk contains $i$ strands. (Trivia:
In his first paper on the subject~\cite{reidemeister}, Reidemeister
put the moves in
another order: He called them {\sl Operation} $1$ ($\Omega_2$), $2$
($\Omega_1$) and $3$ ($\Omega_3$). The terminology
$\{\Omega.1,\Omega.2,\Omega.3\}$ appears in his book
``Knotentheorie''~\cite{reidemeister2} from 1932). 

A more refined classification of Reidemeister moves may be needed
when one wants to prove that a function defined on knot diagrams is a
knot invariant.
There are two common ways in which the
$\{\Omega_1,\Omega_2,\Omega_3\}$-classification can be refined: One is
to take into 
account {\em the orientation of $\rtva$} (that is, to distinguish the
mirror-imaged versions of $\Omega_1$- and $\Omega_3$-moves in
Table~\ref{table:moves}). The other is to also take into account
{\em the orientations of the strands in the changing disk}. 
In this paper a classification is introduced (see
Definition~\ref{def:sameclass} below) which incorporates these
refinements and also the {\em ascendingness} and {\em descendingness}
defined in Section~\ref{sec:ascdesc} below.
The classification in Definition~\ref{def:sameclass} originated
in a study of {\em Gauss diagrams}. In Section~\ref{sec:gauss} the
classification is reformulated in terms of Gauss diagrams.

\subsection{Ascending and descending $\Omega_3$-moves.}\label{sec:ascdesc}
The original part of the classification in
Definition~\ref{def:sameclass} is contained in the
following characterization of an $\Omega_3$-move. We can distinguish the
three strands in the changing disk of an 
$\Omega_3$-move: One is the {\em top} strand, the strand with two
over-crossings. One is the {\em middle} strand, and
one the {\em bottom} strand. When we move along the knot diagram 
in the direction of the orientation, we pass these strands in some
cyclic order.
Call the $\Omega_3$-move {\em descending} if we pass the strands in 
the order {\em top-middle-bottom} and call it {\em ascending} otherwise. 

\begin{theorem}\label{theo:figeight}
The figure eight knot diagram cannot be transformed into its inverse
(the same diagram with the orientation of the strand reversed) without
the use of both ascending and descending $\Omega_3$-moves. In
particular, it cannot be transformed into its inverse without at least
two $\Omega_3$-moves.
\end{theorem}
The proof is given in Section~\ref{sec:indepproofs}.

\subsection{A refined classification of Reidemeister moves.}\label{sec:classif}
By definition, a knot diagram is made from some immersion of $S^1$ in
$\rtva$. We shall assume that the orientation of the knot diagram is
induced from the orientation of $S^1$, so the knot diagram determines
the immersion up to an orientation-preserving diffeomorphism of $S^1$.

Let the knot diagrams $d, k$ be decorated images of the immersions
$i_d,i_k:S^1\to\rtva$. Given these immersions, an equivalence of the
knot diagrams $d$ and $k$
(i.e. an orientation-preserving diffeomorphism of $\rtva$ taking $d$
to $k$) induces an orientation-preserving diffeomorphism of $S^1$.

\begin{deff}[Equivalence of Reidemeister moves]\label{def:sameclass}
Let $X$ be the Reidemeister move that transforms the knot
diagram $d_1$ to $d_2$ inside the changing disk $C$, and let
$Y$ be the Reidemeister move that transforms 
$k_1$ to $k_2$ inside the changing disk $D$. 

$X$ and $Y$ are
{\em directed-equivalent} if
there are {\em orientation-preserving} diffeomorphisms $g,h:C\to D$
such that 
\begin{itemize}
\item[a)]
$g(d_1\cap C)=k_1\cap D$ and $h(d_2\cap C)=k_2\cap D$, 
with the same over-undercrossing information at each crossing.
\item[b)]
Fix some immersions of the oriented circle such
that $d_1,d_2,k_1,k_2$ are their decorated images.
The maps $g,h$ of pieces of knot diagrams induce maps between segments of the
immersing circles. These maps should 
extend to orientation-preserving diffeomorphisms $S^1\to S^1$.
\end{itemize}

$X$ and $Y$ are {\em equivalent} if $X$ is directed-equivalent either to $Y$
or to the reverse move $Y^{-1}$ that changes $k_2$ into $k_1$ inside $D$.
\end{deff}

Every equivalence class splits into exactly two directed equivalence
classes; see Section~\ref{sec:direq}.

Each crossing in a knot diagram is given a {\em sign} by the rule
$\sgn(\ppic{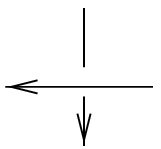})=+1$, $\sgn(\ppic{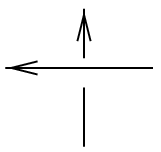})=-1$. 
\begin{deff}[Taxonomy of Reidemeister moves]\label{def:classif}\ 
\begin{itemize}
\item An $\Omega_1$-move is said to be of class $\Omega_{1ij}$,
$i,j=\pm1$. The index $i=+1$
if the kink that the move introduces (or removes) is oriented
counterclockwise, and $i=-1$ otherwise. The index $j$ is the
sign of the crossing that the move introduces or removes.
\item An $\Omega_2$-move is said to be of class $\Omega_{2ij}$,
$i,j=\pm 1$. The index $i=+1$
if the meeting strands are oriented in the same direction (the move is
called {\em direct}), and $i=-1$ otherwise (the move is called {\em inverse}).
The index $j$ is the sign of the affected crossing where the 
overcrossing strand is directed out from the area enclosed by the
meeting strands:
\mbox{$\vcenter{\plainpic{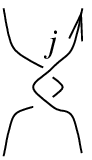}}$} or \mbox{$\vcenter{\plainpic{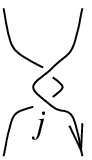}}$}.
\item An $\Omega_3$-move is said to be of class $\Omega_{3ijkm}$,
$i,j,k,m=\pm 1$. The index $m=+1$ if the move is
descending and $m=-1$ otherwise. The indices 
$i,j,k$ are the signs of the crossings in the changing disk:
$$
\vcenter{\plainpic{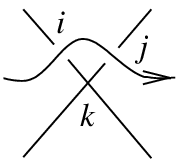}}\rightleftharpoons\vcenter{\plainpic{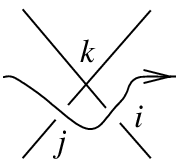}} 
$$
\end{itemize}
\end{deff}
It is easy to verify that Definition~\ref{def:classif} gives the 
splitting of the classes $\Omega_1$, $\Omega_2$ and $\Omega_3$ into 
equivalence classes as defined in Definition~\ref{def:sameclass}.

\subsection{Main results on classes of Reidemeister moves.}\label{sec:statements}
\begin{theorem}\label{theo:smallset}
Any Reidemeister move can be realized by a sequence of moves of the
classes 
$\Omega_{1++}$, $\Omega_{1+-}$, 
$\Omega_{2-+}$, $\Omega_{2--}$, $\Omega_{3-+++}$, $\Omega_{3--+-}$,
followed by a diffeomorphism with support in the changing disk of the
original move.
\end{theorem}
The proof is given in Section~\ref{sec:smallset}. There are 
other sets of six classes of Reidemeister moves that also has this
property, such as the set
$\{\Omega_{1-+}$, $\Omega_{1--}$, $\Omega_{2-+}$, $\Omega_{2+-}$,
$\Omega_{3+--+}$, $\Omega_{3++--}\}$. 

\begin{theorem}\label{theo:indep}
For every knot type there are two knot diagrams, representing knots of
this type, such that to transform them into each other we need at least
\begin{enumerate}
\item one ascending and one descending $\Omega_3$-move, and
\item an $\Omega_{1ij}$-move and an $\Omega_{1kl}$-move,
$(i,j)\ne(k,l)$, $(i,j)\ne (-k,-l)$.
\end{enumerate}
\end{theorem}

Theorem~\ref{theo:indep} extends
Theorem~\ref{theo:figeight} to diagrams of any 
knot type. 
The proof is given in Section~\ref{sec:indepproofs}, and it makes use
of some {\em knot diagram invariants} introduced in Section~\ref{sec:ad}. 
If we remove requirement 1, the remaining statement can be
proved using well-known knot diagram invariants (writhe and winding
number, also cf. Proposition~\ref{prop:regisot}).

\subsection{Organization of the rest of the paper.}
Sections~\ref{sec:direq} and~\ref{sec:ascdescdir} deal with the {\em
direction} of a move.
In Section~\ref{sec:gauss} the definitions and results are reformulated
in terms of {\em Gauss diagrams}. 
Section~\ref{sec:smallset} contains the proof of
Theorem~\ref{theo:smallset}. 

In Section~\ref{sec:ad} several {\em knot diagram
invariants} are introduced, and some of their properties are derived. 
Section~\ref{sec:weirdness} describes how
the simplest of the invariants jumps under a Reidemeister move.
Section~\ref{sec:indepproofs} contains the proofs of
Theorems~\ref{theo:figeight} and~\ref{theo:indep}. The proofs make use
of some of the invariants from Section~\ref{sec:ad}.

\subsection{Directed classes of Reidemeister moves.}\label{sec:direq}
Every equivalence class splits into exactly two directed equivalence 
classes (see Definition~\ref{def:sameclass}). We say that an $\Omega_1$- or
$\Omega_2$-move is in {\em positive direction} if the number of
crossings in the diagram increases. To define the positive direction
of an $\Omega_3$-move we use a construction by V. Arnold
\cite{arnold}:

The three strands in the changing disk (before or after the move) form a
{\em vanishing triangle}. The edges are ordered by the cyclic
order in which we pass them if we move around the diagram in the
direction of the orientation. This cyclic order gives an orientation
of the vanishing triangle. Let $n$ be the number of edges where this
orientation coincides with the  
orientation of the edge. Let $q=(-1)^n$. Then $q=-1$ at one side of the
move and $q=+1$ on the other. We say that the move is in positive
direction if $q=-1$ before the move and $q=+1$ after the move.
\begin{ex}[An $\Omega_3$-move in positive direction]
$$
\ypic{5}{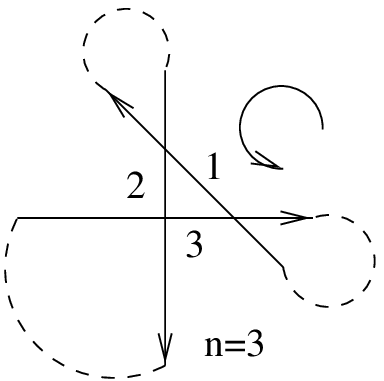}{\LARGE\bf\leadsto}\ypic{5}{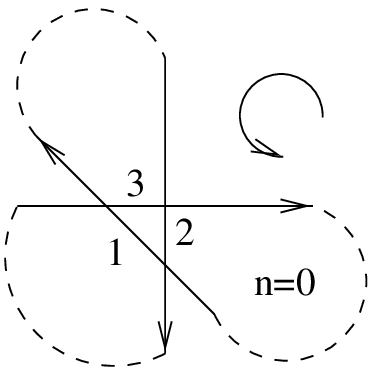}
$$
\end{ex}

\subsection{The direction of ascending and descending
$\Omega_3$-moves.}\label{sec:ascdescdir}
Arnold's definition of the positive direction of a triple point move
appeared in the theory of plane curves. For knot diagrams there is an
ordering of the edges in the vanishing triangle which only
depends on the knot diagram inside the changing disk:
{\em top-middle-bottom}. So it is possible to define a ``positive
direction'' of an $\Omega_3$-move in two different ways: Using
Arnold's ordering of the 
edges, or the top-middle-bottom ordering. The two definitions give
the same ``positive direction'' for descending $\Omega_3$-moves and opposite
``positive direction'' for ascending $\Omega_3$-moves. (In this paper we always
use Arnold's original definition.)

\section{Reidemeister moves of Gauss diagrams.}\label{sec:gauss}
Let the knot diagram $d$ be the decorated image of an immersion of an  
oriented circle, with the orientation of the diagram given by the
orientation of the circle. 
The {\em Gauss diagram} $G_d$ of $d$ is the immersing oriented  
circle with the preimages of each crossing point connected by a signed
arrow. The arrow points from overpass to underpass and is given the
sign of the crossing.
Gauss diagrams are considered up to orientation-preserving
diffeomorphism of the circle, so equivalent knot diagrams have the
same Gauss diagram.  
\begin{ex}[The Gauss diagram of a knot diagram.]
$$
\ypic{4}{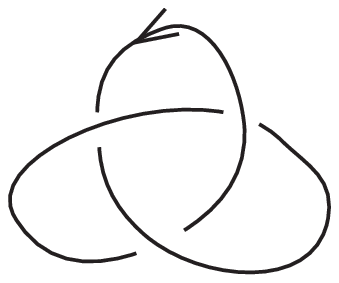}\ \ \ \ypic{4}{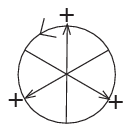}
$$
\end{ex}
The Gauss diagram determines a knot diagram on the sphere
up to an orientation-preserving diffeomorphism of the sphere.
That is, the Gauss diagram determines a knot diagram in the plane up to a 
choice of a point of infinity outside the diagram on the sphere and up to
equivalence of knot diagrams. 

A Reidemeister move induces a transformation of Gauss diagrams. 
An $\Omega_i$-move affects $i$ arrows, with end points on $i$ segments of the
circle of the Gauss diagram, without changing the rest of the diagram.

The original motivation for the choice of the definitions in this paper was
the following obvious Proposition:
\begin{prop}
The Reidemeister moves $X:d_1\mapsto d_2$ and $Y:k_1\mapsto k_2$ are
directed-equivalent if and only if the 
abstract Gauss diagrams obtained by erasing all arrows in
$G_{d_1}$ and $G_{k_1}$ that are not affected by the move are the same
(up to an orientation-preserving diffeomorphism of the circle), and
likewise for $G_{d_2}$ and $G_{k_2}$.\hfill$\Box$
\end{prop}
That is, $X$ and $Y$ are directed-equivalent if and only if
the induced transformations of the Gauss diagrams look the same when
disregarding the unchanged part of the diagram.

\begin{cor}
A function on Gauss diagrams defines a knot invariant if and only if it is
invariant under the transformations in Table~\ref{table:gaussmoves},
where the dotted
segments indicate a part of the Gauss diagram that is unchanged.
\end{cor}
In Table~\ref{table:gaussmoves}, the positive direction of a move is
from left to right.
\begin{table}\caption{A sufficient set of Reidemeister moves for Gauss diagrams.}\label{table:gaussmoves}
\def\epsfsize#1#2{0.8#1} 
$$
\begin{array}{cccccc}
\plainpic{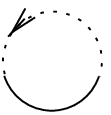}&\rightleftharpoons\atop\Omega_{1++}&\plainpic{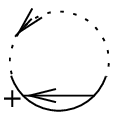}&
\plainpic{o1xx-.eps}&\rightleftharpoons\atop\Omega_{1+-}&\plainpic{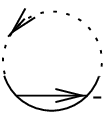}\\
\plainpic{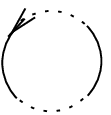}&\rightleftharpoons\atop\Omega_{2-+}&\plainpic{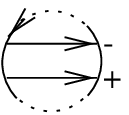}&
\plainpic{o2xx-.eps}&\rightleftharpoons\atop\Omega_{2--}&\plainpic{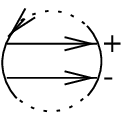}\\
\plainpic{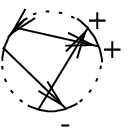}&\rightleftharpoons\atop\Omega_{3-+++}&\plainpic{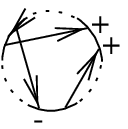}&
\plainpic{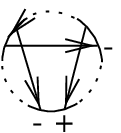}&\rightleftharpoons\atop\Omega_{3--+-}&\plainpic{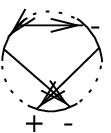}
\end{array}
$$
\def\epsfsize#1#2{\ifnum\value{pictype}=0 0.7#1\else\epsfxsize\fi} 
\end{table}

\section{Proof of Theorem~\ref{theo:smallset}.}\label{sec:smallset}
The proof of Theorem~\ref{theo:smallset} consists of three lemmas. The
method is to draw movies that show explicitly how a
Reidemeister 
move of a class we want to exclude is replaced by a sequence of other
Reidemeister moves. All moves in this sequence take place inside
the changing disk of the move we replace.

\begin{lemma}\label{lemma:om1repl}
Any $\Omega_1$-move can be replaced by
$\Omega_{1++}$-, $\Omega_{1+-}$- and $\Omega_{2--}$-moves.
\end{lemma}
\subsubsection*{Proof:}
We can replace an $\Omega_{1-+}$-move with the sequence
$$
\ypic{3}{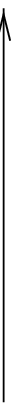}\ {\rightharpoonup\atop \Omega_{2--}}\ 
\ypic{3}{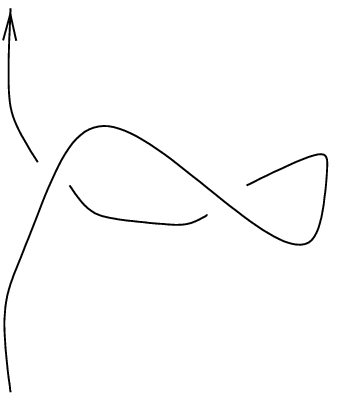}\ {\rightharpoonup\atop \Omega_{1+-}}\ 
\ypic{3}{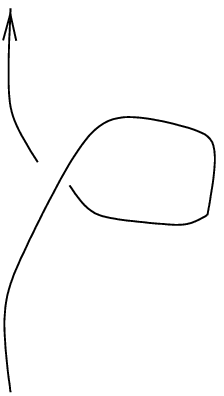}.
$$
And we can replace an $\Omega_{1--}$-move with the sequence
$$
\ypic{3}{psi1seq-+1.eps}\ {\rightharpoonup\atop \Omega_{2--}}\ 
\ypic{3}{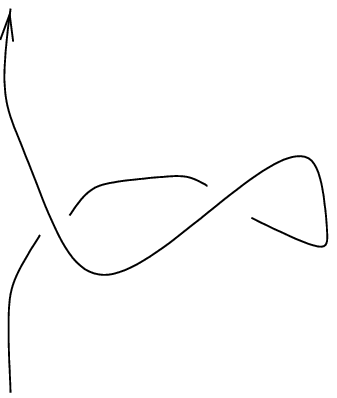}\ {\rightharpoonup\atop \Omega_{1++}}\ 
\ypic{3}{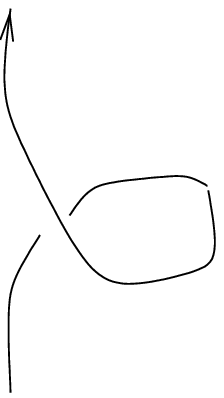}.
$$
\hfill$\Box$

\begin{lemma}\label{lemma:om2repl}
Any direct $\Omega_2$-move can be replaced by
$\Omega_{2--}$-, $\Omega_{1-+}$-, $\Omega_{1--}$-, 
$\Omega_{3-+++}$- and $\Omega_{3--+-}$-moves.
\end{lemma}
\subsubsection*{Proof:}
We can replace an $\Omega_{2++}$-move with the sequence
$$
\ypic{3}{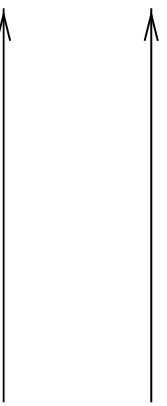}\ {\rightharpoonup\atop \Omega_{1-+}}\ 
\ypic{3}{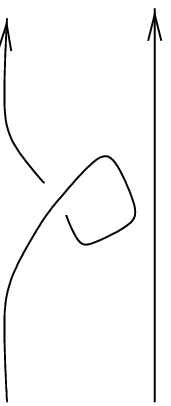}\ {\rightharpoonup\atop \Omega_{2--}}\ 
\ypic{3}{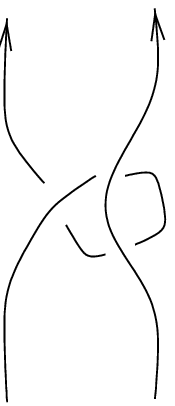}\ {\rightharpoonup\atop \Omega_{3-+++}}\ 
\ypic{3}{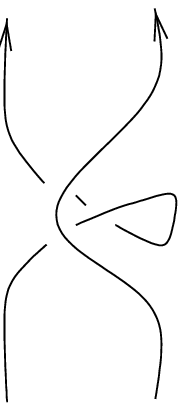}\ {\rightharpoonup\atop \Omega_{1-+}}\ 
\ypic{3}{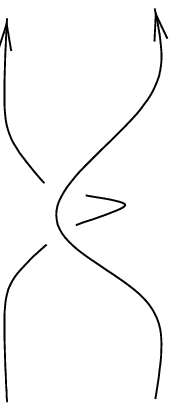}.
$$
In a similar way we can replace an $\Omega_{2+-}$-move with a sequence
of moves of class $\Omega_{1--}$, $\Omega_{2--}$, $\Omega_{3--+-}$ and
$\Omega_{1--}$.\hfill$\Box$

\begin{lemma}\label{lemma:om3repl}
Any ascending (descending) $\Omega_3$-move can be
replaced by $\Omega_2$-moves and an $\Omega_3$-move of any ascending
(descending) class.
\end{lemma}
\subsubsection*{Proof:}
We can replace an $\Omega_3$-move with a sequence which looks like
$$
\ypic{4}{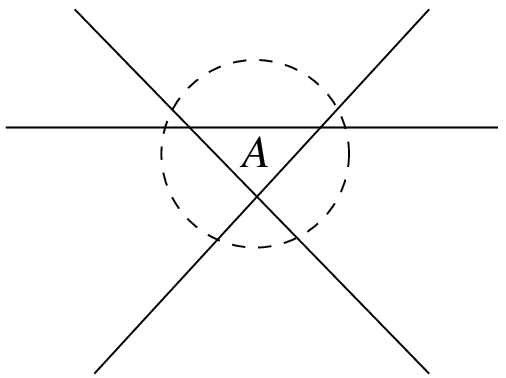}\ {\rightharpoonup\atop \Omega_{2}}\ 
\ypic{4}{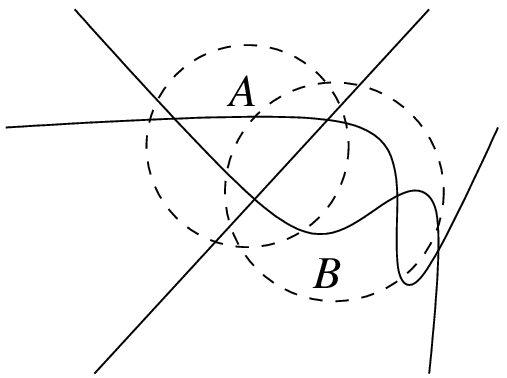}\ {\rightharpoonup\atop \Omega_{3}}\ 
\ypic{4}{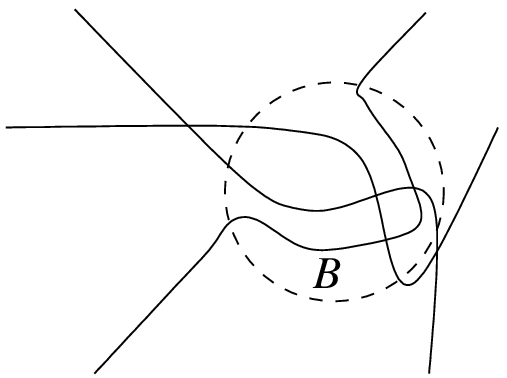}\ {\rightharpoonup\atop \Omega_{2}}\ 
\ypic{4}{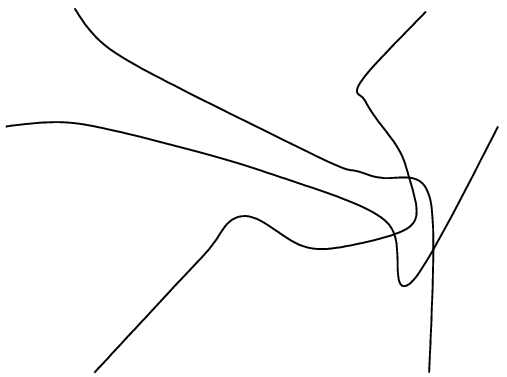}.
$$
if we temporarily forget about orientation, 
over/under-crossing information and the ordering of the edges.
Here we replace the $\Omega_3$-move in the changing disk $A$ with
a move in the changing disk $B$. 
The $\Omega_3$-moves in $A$ and in $B$ are either both ascending or both
descending. The first three signs $i,j,k$ in the
$\Omega_{3ijkm}$ index are given by the signs of the three crossings
in the changing disk as indicated in Definition~\ref{def:classif}.

We have three choices of wedges where to make the first $\Omega_2$-move. If we make
the move between the middle and bottom strands, $i$ and $j$ change
place in the index 
and $k$ is reversed: $(i,j,k)\mapsto (j,i,-k)$. If we make the move between 
the top strand and the strand that
crosses the top strand at the $i$-crossing then 
$(i,j,k)\mapsto (j,-i,k)$. If we make the move between the top strand
and the strand that crosses the top strand at the $j$-crossing we get
$(i,j,k)\mapsto (-j,i,k)$. Repeating this process with the new changing
disk, we can get any index. This completes the proof of
Theorem~\ref{theo:smallset}.\hfill$\Box$

\section{Knot diagram invariants.}\label{sec:ad}
\subsection{Writhe and winding number.}
Note that by a {\em knot diagram invariant} we mean a function on knot
diagrams that 
is unchanged under equivalence of knot diagrams. There are some
well-known knot diagram invariants which are not knot invariants, 
foremost the {\em writhe} and the {\em
winding number}. The writhe is the sum of the signs of all 
crossings in the diagram. The winding number is the degree of the map
$S^1 \to S^1$ taking a point on the immersing circle to the direction
of the tangent vector in the plane. 

The main previously known result about classes of Reidemeister moves
is Proposition~\ref{prop:regisot} below. Note that the writhe and the
winding number are unchanged under $\Omega_1$- and $\Omega_2$-moves,
which proves the easy ``only if''-part of the Proposition.

\begin{prop}[B. Trace~\cite{trace}]\label{prop:regisot}
Two knot diagrams of isotopic knots can be transformed into each other
without the use of $\Omega_1$-moves (they are {\em regularly isotopic})
if and only if they have the same writhe and winding number.
\end{prop}

\subsection{Connected sum of knot diagrams.}\label{sec:consum}
The knot diagram $k$ is a {\em connected sum} of $k_1$ and $k_2$ if 
$k$ can be split by a simple closed curve intersecting the diagram in
two points, into diagrams equivalent to $k_1$ and $k_2$, as below:
$$
\ypic{10}{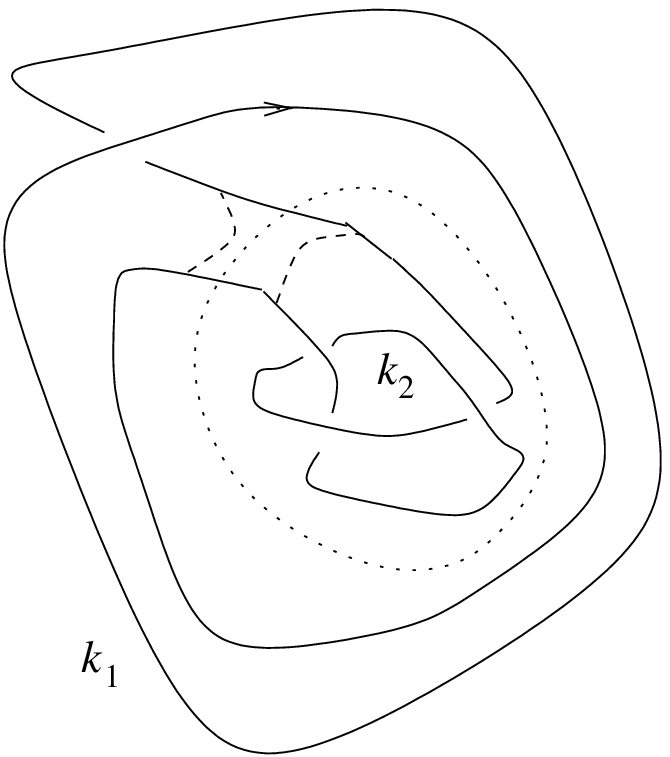}
$$
If $k$ is a connected sum of $k_1$ and $k_2$, the corresponding Gauss
diagram $G_k$ is a connected sum of $G_{k_1}$ and $G_{k_2}$ in the
sense that $G_k$ can be split into $G_{k_1}$ and $G_{k_2}$: In the
example above, the terms have Gauss diagrams
$G_{k_1}=\ppic{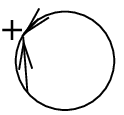}$ and $G_{k_2}=
\ppic{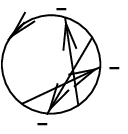}$, while the sum has Gauss diagram \ppic{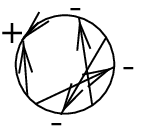}.
A knot diagram invariant $V$ is called
{\em additive} if $V(k)=V(k_1)+V(k_2)$ whenever $k$ is a connected
sum of $k_1$ and $k_2$. 
\begin{lemma}\label{lemma:localinvproperties}
The writhe is additive. The winding number is not, but
if $k$ is a connected sum of $k_1$ and $k_2$, then 
$$
\mbox{winding number}(k)=\mbox{winding number}(k_1)+\mbox{winding
number}(k_2)+i
$$
where $i$ depends on how the sum is done: $i=+1$ if we cut like
\cpic{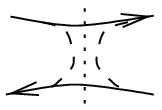} and $i=-1$ if we cut like
\cpic{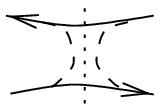}. 
\end{lemma}
\subsubsection*{Proof:}
The first part is evident. For the second part: If we cut
the knot diagram $k$ into $k_1$ and $k_2$ we introduce two new
half-turns of the direction vector, in the same direction. Hence the
sum of the winding numbers 
of $k_1$ and $k_2$ differ by $+1$ from the winding number of their connected
sum $k$ if the half-turns are in positive direction, and by $-1$ if they
are in negative direction.\hfill$\Box$

\subsection{The knot diagram invariants $A_n$, $D_n$ and $W_n$.}
For a knot diagram $k$ let $A_4(k)$ be the sum of the {\em signs} of
all subdiagrams of the form $a_4=\ppic{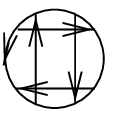}$ in the
Gauss diagram $G_k$ of $k$. A subdiagram is an abstract Gauss diagram
that can be created from $G_k$ by removing arrows.
The sign of the subdiagram is the product
of the signs of the arrows in the subdiagram. In the same way define
$D_4$ as the sum of the signs of the subdiagrams of the form
$d_4=\ppic{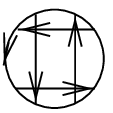}$. In Section~\ref{sec:indepproofs} we use $A_4$
and $D_4$ to prove Theorems~\ref{theo:figeight} and~\ref{theo:indep}. 

More generally, we
define $A_n(k)$ and $D_n(k)$, $n\ge 4$, as the sum of the signs of the
subdiagrams of $G_k$ with $n$ arrows arranged cyclically as
$$\xpic{0.5}{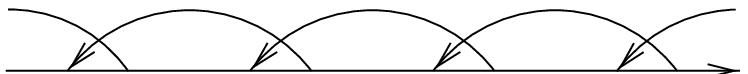}$$ respectively $$\xpic{0.5}{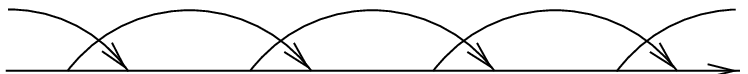}.$$
Define $W_n(k)$, $n=3,5,7,\ldots$, as the sum
of the signs of the subdiagrams of $G_k$ with $n$ arrows which are diameters
directed so that no two arrow heads are adjacent on the circle. That
is, $W_5$ is the sum of the signs of the subdiagrams that look like
$w_5=\ppic{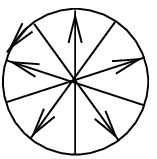}$. 

Obviously, $A_n$, $D_n$ and $W_n$ are knot diagram invariants.
The abstract subdiagrams that define these invariants, which we call
$a_n$, $d_n$ and $w_n$, are examples of
{\em arrow 
diagrams}. Functions defined in this way were introduced by M. Polyak
and O. Viro~\cite{viro}.

\begin{prop}[Properties of $A_n$, $D_n$ and $W_n$]\label{prop:invarproperties}
$A_n$, $D_n$ and $W_n$ are
additive under the connected sum of knot diagrams. $A_n(k)$ and $D_n(k)$ are
multiplied by $(-1)^n$ by mirroring of $k$ (change of orientation of $\rtva$)
and change into each other under inversion of $k$ (change of orientation of
the strand). $W_n(k)$ is multiplied by $-1$ by mirroring of $k$ and is
unchanged by inversion of $k$.
\end{prop}
\subsubsection*{Proof:}
Mirroring of the knot diagram affects the Gauss diagram by reversing
the signs of all arrows. The sign of an $a_n$-, $d_n$- or
$w_n$-subdiagram is a product of the signs of $n$ arrows, so every
sign is multiplied by $(-1)^n$. For a $w_n$-subdiagram, $n$ is always odd.
 
Inversion of the knot diagram affects the Gauss diagram by reversing
the orientation of the circle, so the statement is evident.

Additivity is because the Gauss
diagram of a connected sum is a connected sum of Gauss diagrams (see
Section~\ref{sec:consum}), and each $a_n$-, $d_n$-, or
$w_n$-subdiagram evidently must belong to one of the summands.\hfill$\Box$ 

\begin{prop}[Reidemeister move invariance of $A_n$, $D_n$ and $W_n$]\label{prop:invarinvar}
$A_n$, $D_n$ and $W_n$ do not
change under $\Omega_1$- and $\Omega_2$-moves. $A_n$ does not change
under descending $\Omega_3$-moves and $D_n$ does not change under
ascending $\Omega_3$-moves.
\end{prop}
We assert that $A_n$ ($D_n$) does change under some ascending
(descending) $\Omega_3$-moves and $W_n$ does change under some
$\Omega_3$-moves of both ascending and descending type.
(That is, they are not knot invariants.) 

\subsubsection*{Proof of $\Omega_1$-invariance:}
An $\Omega_1$-move introduces (or removes) one arrow in the Gauss
diagram as in Table~\ref{table:gaussmoves}. This arrow is {\em isolated};
its head and tail are adjacent on the 
circle.  $a_n$-, $d_n$- and $w_n$-subdiagrams contain no such arrows,
so the set of such subdiagrams does not change.\hfill$\Box$
\subsubsection*{Proof of $\Omega_2$-invariance:}
An $\Omega_2$-move introduces (or removes)
two arrows in the Gauss diagram as in
Table~\ref{table:gaussmoves}. The tails of these arrows are 
adjacent on the circle, and the heads are also adjacent on the
circle. This is not the case for any two arrows in an $a_n$-, $d_n$- or
$w_n$-subdiagram. There may be $a_n$-, $d_n$- or $w_n$-subdiagrams
that contain one 
of the two arrows affected by the move. However,
such subdiagrams come in pairs with opposite sign, since the two arrows have
opposite sign and point from adjacent points on the circle to
adjacent points on the circle. Hence $A_n$, $D_n$ and $W_n$ are
invariant under $\Omega_1$- and $\Omega_2$-moves.\hfill$\Box$
\subsubsection*{Proof of invariance of $A_n$ ($D_n$) under descending
(ascending) $\Omega_3$-moves:} 
By Lemma~\ref{lemma:om3repl}, a knot diagram invariant that is
unchanged under $\Omega_1$- and $\Omega_2$-moves is invariant under
all ascending (respectively descending) $\Omega_3$-moves if it is
invariant under moves of class $\Omega_{3--+-}$ ($\Omega_{3-+++}$).
The corresponding moves of Gauss diagrams are depicted in
Table~\ref{table:gaussmoves}. 

An $\Omega_3$-move affects three arrows in the Gauss diagram.
The move obviously does not change the number of $a_n$- and
$d_n$-subdiagrams that only 
contain one of these three arrows. The subdiagrams that contain two of
the affected 
arrows might change. These two arrows must then have adjacent end
points at some point on the circle. 

Any two arrows in an $a_n$-subdiagram have
adjacent end points at no more than
one place on the circle. Then these end points are one head and one
tail. If the two arrows are crossed, 
the tail comes before the head on the circle, and if the arrows are
not crossed, the 
head comes before the tail. 

We see in Table~\ref{table:gaussmoves} that
for an $\Omega_{3-+++}$-move, two affected arrows that have
a head and a tail next to each other are crossed if and only if the head comes
before the tail. Hence no pair of
arrows that are 
affected by an $\Omega_{3-+++}$-move can belong to an $a_n$-subdiagram.
Consequently $A_n$ is invariant under descending
$\Omega_3$-moves. In the same way we see that no pair of arrows that
are affected by an $\Omega_{3--+-}$-move can belong to a
$d_n$-subdiagram.\hfill$\Box$

\subsection{The Weirdness $W_3$. A geometric formula for the jump.}\label{sec:weirdness}
The function $W_3$, which is given by subdiagrams on the form
$w_3=\ppic{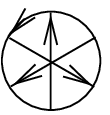}$, 
is called {\em Weirdness}. We shall give a formula for the
jump of $W_3$ under a positively directed $\Omega_3$-move.

Let $k$ be a knot diagram and $C$ a selected changing disk for a
$\Omega_3$-move on $k$. The piece of knot diagram exterior to $C$
consists of three strands, which we label $a$, $b$ and $c$. $a$ is the
strand that connects the top strand inside $C$ with the middle strand.
$b$ connects the middle and the bottom strand, and $c$ connects the
top and the bottom strand.
By the {\em local linking number $llk$} of the pair $(k,C)$
we mean the following (cf. Example~\ref{ex:llk} below):

Erase the piece of knot diagram that lies inside $C$. Erase the
$c$ strand. Make the remaining $a$ and $b$ strands into a
two-component link diagram by connecting the two ends of each strand
with a strand inside $C$ in such a way that the strand that closes the
$b$ strand {\em does not pass over} the strand that closes the $a$ strand.

The local linking number is the linking number of this link. (The
linking number is an invariant of two-component links. If we pick one
of the components, the linking number is the sum of the signs of the
crossings where this component crosses over the other.)

$llk$ is clearly invariant under equivalence of knot diagrams (as
long as we keep track of the selected changing disk.)
\begin{prop}\label{prop:llk}
$W_3$ jumps by $llk(k,C)$ under a positively directed $\Omega_3$-move
in the changing disk $C$ on the knot diagram $k$.
\end{prop}

\begin{ex}[Computing the local linking number.]\label{ex:llk}\rm
Consider the changing disk $C$ in the knot diagram
$$
k=\ypic{4}{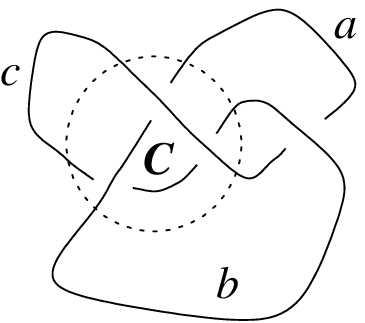}.
$$
We make the knot diagram into a link diagram as below:
$$
\ypic{4}{llk1.eps}\ \longmapsto\ \ypic{4}{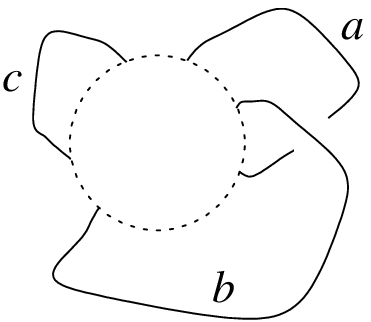}\ \longmapsto\ \ypic{4}{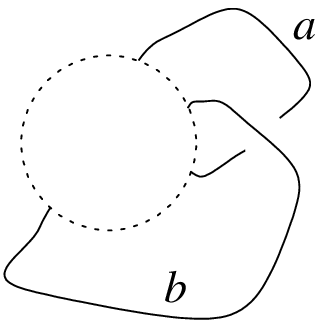}\ \longmapsto\ \ypic{4}{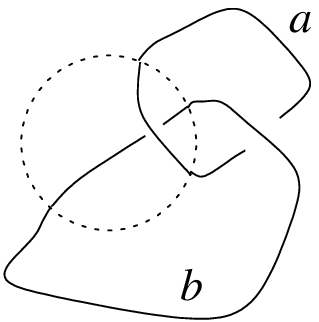}
$$
$llk$ is the linking number of this link.
\end{ex}

\subsubsection*{Proof of Proposition~\ref{prop:llk}:}
By Lemma~\ref{lemma:om3repl}, we can replace an $\Omega_3$-move with a
sequence consisting of $\Omega_2$-moves and a move of any
$\Omega_3$-class with the same ascendingness/descendingness. The main
trick is to make a move in the changing disk $B$ instead of the
changing disk $A$ in the picture below:
$$
k_1=\ypic{4}{om3redseq1.eps}\ {\rightharpoonup\atop \Omega_{2}}\ 
\ypic{4}{om3redseq2.eps}=k_2
$$
Here $llk(k_1,A)=llk(k_2,A)$ since $\Omega_2$-moves outside $A$ cannot
change the linking number. $llk(k_2,A)=llk(k_2,B)$ since the most that
can happen is that an overcrossing of the $a$ (or $b$) strand is
replaced with an overcrossing of the $b$ ($a$) strand with opposite
sign. This does not change the linking number.
Hence it is sufficient to prove the statement for
$\Omega_{3-+++}$-moves and $\Omega_{3--+-}$-moves. 
We consider
$w_3$-subdiagrams of a Gauss diagram, such that two of the arrows in
the subdiagrams belong to the changing part. (Subdiagrams such that
only one arrow belong to the changing part are unchanged by the move.)
The changing part is
depicted in Table~\ref{table:gaussmoves}. We see that there can be no
subdiagrams of this kind after the positively directed move, when the
changing part looks as on the right side in the Table. Before the
positively directed $\Omega_{3-+++}$-move we can have 
$w_3$-subdiagrams that look like \\ \ppic{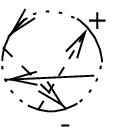}, which contains the
two dotted arrows in the changing part and one arrow pointing between
unchanged segments as shown. Before the
positively directed $\Omega_{3--+-}$-move we can have 
$w_3$-subdiagrams that look like \ppic{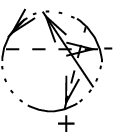}. 

In both cases, the sign of the subdiagram is minus the sign of the 
third arrow in the subdiagram. The third arrow points from the segment
that represents 
the $b$ strand to the segment that
represents the $a$ strand. So when we make the move in positive
direction, we add to $W_3$ the sum of the signs of the crossings where
the $b$ strand crosses over the $a$ strand. This is the
local linking number, since we close the link without introducing any
new overcrossings of $b$ over $a$.\hfill$\Box$

\section{Proof of Theorems~\ref{theo:figeight} and
\ref{theo:indep}.}\label{sec:indepproofs} 

\subsection{Proof of Theorem~\ref{theo:figeight}.}
We shall prove that the figure eight diagram $d_8$ in
Figure~\ref{fig:figeight} cannot be transformed into its inverse
without the use of both ascending and descending $\Omega_3$-moves. 
$d_8$ has Gauss diagram $G_{d_8}=\ppic{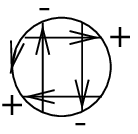}$.
Obviously $A(d_8)=1$ while $D(d_8)=0$. The inverse knot diagram has 
$A=0$ and $D=1$ by Proposition~\ref{prop:invarproperties}. The result
follows since $A$ jumps only under 
ascending $\Omega_3$-moves and $D$ jumps only under descending
$\Omega_3$-moves.\hfill$\Box$

\subsection{Proof of Theorem~\ref{theo:indep}.}
For each knot diagram, we shall construct a knot diagram of an
isotopic knot, such that the two knot diagrams can only be transformed
into each other with the use of both descending and ascending
$\Omega_3$-moves, and move of some classes $\Omega_{1ij}$ and
$\Omega_{1kl}$, $(i,j)\ne\pm (k,l)$.
Consider the knot diagram $a$ of the unknot:
$$
{\Huge a=}\ypic{5}{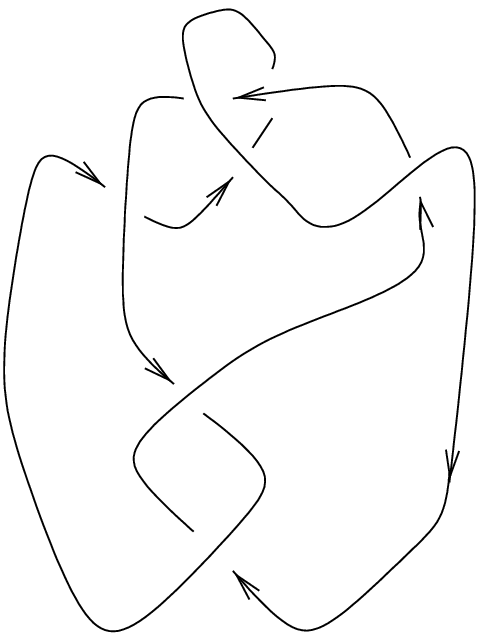} 
$$

$a$ has Gauss diagrams $G_a=\ypic{3}{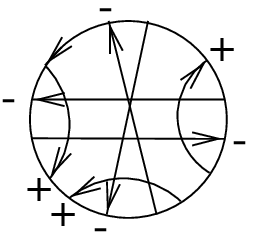}$. We see that
$D(a)=-1$, $A(a)=0$, $\mbox{writhe}(a)=-1$ and $\mbox{winding
number}(a)=0$. The inverse knot diagram $\tilde{a}$ has
$D(\tilde{a})=0$, $A(\tilde{a})=-1$, $\mbox{writhe}(a)=-1$ and
$\mbox{winding number}(a)=0$. (The writhe is unchanged under
inversion, and the winding number is multiplied with $-1$.) For every
knot diagram $d$ there is some
knot diagram $k$ that is a connected sum of $d$, $a$ and
$\tilde{a}$. This connected sum should be made at corresponding points
on $a$ and $\tilde{a}$.

$k$ and $d$ describe isotopic knots, since we have
just added two unknots.
$A(k)=A(d)-1$ and $D(k)=D(d)-1$ by
Proposition~\ref{prop:invarproperties}, hence we need both ascending
and 
ascending $\Omega_3$-moves to transform $k$ to $d$.
$\mbox{writhe}(k)=\mbox{writhe}(d)-2$ and
$\mbox{winding number}(k)=\mbox{winding number}(d)$ by
Lemma~\ref{lemma:localinvproperties} (because the
contribution in winding number from $a$ and $\tilde{a}$ cancel).
The winding number and writhe jump by $\pm i$ and $\pm j$ under an
$\Omega_{1ij}$-move and are invariant under $\Omega_2$- and
$\Omega_3$-moves, so a collection of moves of classes $\Omega_{1ij}$ and
$\Omega_{1,-i,-j}$ can never raise the writhe by $2$ and leave the
winding number unchanged.\hfill$\Box$

\end{document}